\documentclass[leqno]{amsart}

\usepackage{calc}
\usepackage{geometry}  
\usepackage{graphicx}
\usepackage{amsmath,amssymb,epsfig,mathtools}
\usepackage{amsthm}
\usepackage{enumerate}
\usepackage{caption}
\usepackage{tikz}
\usepackage{hyperref}

\makeatletter
\renewcommand*\env@matrix[1][*\c@MaxMatrixCols c]{%
  \hskip -\arraycolsep
  \let\@ifnextchar\new@ifnextchar
  \array{#1}}
\makeatother

\newcommand{\PP}{{\mathbf P}}

\newcommand{\ZZ}{{\mathbb Z}}

\newcommand{\E}{{\mathcal E}}

\newcommand{\Max}{{ \mathrm{Max} }}

\newcommand{\Prob}{{ \mathrm{Prob} }}

\newcommand{\FF}{\mathbb{F}}

\newcommand{\QQ}{\mathbb{Q}}

\newtheorem{theorem}{\bf{Theorem}}[section]

\theoremstyle{definition}

\newtheorem{remark}[theorem]{\bf{Remark}}

\def\MoW{Mordell--\mbox{\kern-.12em}Weil}
\def\0{^{\phantom0}}
\def\NS{{\mathop{\rm NS}\nolimits}}
\def\disc{{\mathop{\rm disc}\nolimits}}

\begin{document}

\title{New rank records for elliptic curves having rational torsion}
\author{Noam D. Elkies}
\address[N.~Elkies]{Dept. of Mathematics, Harvard University, 1 Oxford St., Cambridge, MA 02138}
\email{elkies@math.harvard.edu}
\thanks{Elkies was supported by NSF grants DMS-0501029, DMS-1100511, and DMS-1502161,
a Radcliffe Fellowship, and the Simons Collaboration on
Arithmetic Geometry, Number Theory, and Computation.}

\author{Zev Klagsbrun}
\address[Z.~Klagsbrun]{Center for Communications Research, 4320 Westerra Court, San Diego, CA 92121}
\email{zdklags@ccrwest.org}

\bibliographystyle{alpha}

\begin{abstract}
We present rank-record breaking elliptic curves having torsion subgroups $\ZZ/2\ZZ$, $\ZZ/3\ZZ$, $\ZZ/4\ZZ$, $\ZZ/6\ZZ$, and $\ZZ/7\ZZ$.
\end{abstract}
\maketitle

\section{Introduction}
\label{sec:intro}

Given an elliptic curve $E/\QQ$, the \MoW\ theorem states that
the group of rational points $E(\QQ)$ is isormorpic to $\ZZ^r \times T$,
where $r$ is the \emph{rank of $E$} and $T$ is a finite group called
the \emph{torsion subgroup} of $E$~\cite{MWThm}.
While the groups that can appear as $T$ have been fully characterized by
Mazur~\cite{Mazur}, which ranks occur is a question that goes back to
Poincar\'{e}~\cite{Poincare} and has been the subject of competing
folklore conjectures.

One side, claiming that ranks are bounded, has recently been bolstered by
several different models \cite{WatkinsTwist, Watkins, PPVW}
that predict that all but finitely many elliptic curves have rank
at most $21$, with stronger conjectured bounds on which ranks occur
infinitely often for each possible torsion group~$T$.
(For example, if $T=\ZZ/n\ZZ$ for $n=2,3,\ldots,8$ then the bound $21$
is replaced by $13,9,7,5,5,3,3$.)
The other side, arguing that ranks are unbounded,
has relied on periodically exhibiting curves of larger and larger rank.

Our work continues that tradition, exhibiting rank-record breaking curves
for the torsion subgroups $\ZZ/2\ZZ$, $\ZZ/3\ZZ$, $\ZZ/4\ZZ$, $\ZZ/6\ZZ$,
and $\ZZ/7\ZZ$, which constitute one-third of the $15$ groups that Mazur
showed can appear as the torsion subgroup of an elliptic curve over $\QQ$.

At the same time, our work provides, at best, limited evidence that
ranks are unbounded.  We broke five different records, and found
numerous new curves whose ranks tie the old records
(and many more whose ranks exceed the heuristically conjectured
asymptotic upper bounds).  But the scale of this search was
vastly larger than any previously attempted, and yet
we could not break any of the previous records by more than~$1$,
and in each case found only a handful of curves
(in most cases a single curve) with the new record rank. This suggests that
the growth of ranks of elliptic curves might indeed peter out at some point.

\subsection{Organization}
This paper largely splits into three parts.
The first consists of Sections \ref{sec:MestreNagao}~--~\ref{sec:choosingBC},
which describe the methods that we used to search for curves of high rank,
as well as Section~\ref{sec:openQ}, which presents some open questions about
our methods. The second, Sections \ref{sec:MainResults} -- \ref{sec:Z7Z},
describes our results, including details of our searches in each of
the torsion families considered. Section \ref{sec:Z2Z} also includes
a previously unpublished family of elliptic K3 surfaces $\E_u/\QQ(t)$
that have \MoW\ group $\ZZ/2\ZZ \times \ZZ^9$ for each $u \ne \pm 1, \pm 2$
for which $5-u^2$ is a square.
We exhibit generators for $\E_u(\QQ(t))$ in Appendix~\ref{app:pointsonShimura}.
The third and final part of this paper consists of
Appendix~\ref{app:modelsandpoints}, which presents models for the record-breaking
curves we discovered and points that generate their \MoW\ groups.

\section{The method of Mestre and Nagao}
\label{sec:MestreNagao}

The core ingredient in our search was a well-known method,
originally conceived by Mestre,
for finding elliptic curves having large \MoW\ rank. 
We start with an elliptic fibration
$\E/\QQ(t)$ having \MoW\ rank $r$,
and then attempt to find good values of $t$ for which
the specialization $E_t$ has particularly large rank~\cite{Mestre15}.

A theorem of Silverman \cite{SilvermanSpecialization} states that
all but finitely many specializations $E_t$ of $\E$ have rank at least $r$,
so this approach effectively gives us $r$ independent rational points on
each specialization for free. 

The method for finding values of $t$ for which the rank of $E_t$
is significantly larger than $r$ has its roots in the observation of
Birch and Swinnerton-Dyer that curves that have unusually many points
modulo $p$ for most $p$ should have many rational points as well~\cite{BSD}
and in Mestre's work on Weil's explicit formula for
elliptic curves~\cite{MestreExplicit}.
The idea is to construct a score $S(t,B)$ that incorporates
the number of points $N_p(E_t)$ on $\overline{E_t}(\FF_p)$
for all primes $p \le B$ where $E_t$ has good reduction,
and then to search for rational points on $E_t$ for those values of $t$
in a search region for which $S(t,B)$ is above some threshold.
While this basic method was first used by Mestre to find the first
curves over $\QQ$ having rank $12$~\cite{Mestre12Q},
its first use in a family $\E/\QQ(t)$ appears to be due to Nagao~\cite{Nagao17}.

Nagao considered the scores 
$$
S_1(t,B) = \sum_{
     \substack{p < B \\ E_t \text{ has good} \\ \text{reduction  at } p}
   }
   \frac{-a_p(E_t) + 2}{N_p(E_t)}\log p
\text{ and }
S_2(t,B) = \frac{1}{B}\sum_{
     \substack{p < B \\ E_t \text{ has good} \\ \text{reduction  at } p}
   }
   -a_p(E_t) \log p,
$$
which, when large, suggest via Weil's explicit formula for elliptic curves
\cite{MestreExplicit} that the order of the vanishing of the $L$-function
$L_{E_t}(s)$ at $s = 1$ should be large as well.

We choose to evaluate a different sum
\begin{equation}
\label{eq:partsum}
S(t,B) = \sum_{
     \substack{p < B \\ E_t \text{ has good} \\ \text{reduction  at } p}
  } \log\left ( \frac{N_p(E_t)}{p} \right)
\end{equation}
as in \cite{Elkies3Lectures}, so that $\exp(-S(t,B))$
is a partial product for 
\begin{equation}
\label{eq:partprod}
\prod_{ \substack{p < B \\ E_t \text{ has good} \\ \text{reduction  at } p}}
  \!\! (1 - a_p(E_t)p^{-s}  + p^{1-2s})^{-1}
\end{equation}
of the Euler product for $L_{E_t}(s)$ evaluated at $s = 1$
(ignoring the finitely many factors at primes of bad reduction).
The conjecture of Birch and Swinnerton-Dyer
suggests that when $E_t$ has large rank
such partial products should rapidly approach zero,
and thus that $S(t,B)$ should be large.

\section{Computational Techniques}

Computing any of the sums in Section \ref{sec:MestreNagao} 
would be computationally infeasible
for a large range of $t$ if one needed to individually compute
$a_p(E_t)$ for each $p < B$ and each value of $t$.
To scale Mestre's method to extremely large search regions,
we took advantage of three computational tricks.

First, as observed by Nagao \cite{Nagao},
$a_p(t)$ depends only on $t \pmod{p}$.
As a result, one can first compute $a_p(t)$
for all $p \le B$ and for all $t \in \FF_p$ for which $\Delta_{E_t} \neq 0$,
and then use the pre-computed values to calculate
$S(t,B)$ for each $t$ in the search region.

The second trick, also due to Nagao \cite{Nagao},
lets us concentrate our computation on the most promising values of $t$.
Rather than compute $S(t,B)$ for all $t$ in the search region,
we choose an increasing series of bounds $B_0 \le B_1 \le \ldots \le B_m = B$
and cutoffs $C_0 \le C_1 \le \ldots \le C_m = C$,
and only compute $S(t,B_i)$ for $i \ge 1$
for those values of $t$ for which $S(t,B_j) \ge C_j$ for all $0 \le j < i$.

These first two tricks appear to be well-known
(see \cite{Fermigier}, for example).
The third trick, which is apparently due to Elkies \cite{Elkies3Lectures},
seems to be less widely known, and we describe it in detail below.

\subsection{Sieving}

Rather than computing $S(t,B)$ for each value of $t$ by looking up
the values of $N_p(t)$ (or more likely, $\log (N_p(t)/p)$)
for each prime $p < B$, sieving computes $S(t,B)$ for a large number 
of values of $t = a/b$ at once. The algorithm works as follows:

Fix a value of $b$ and an interval $[a_0, a_0 + N)$.
We allocate a counter array $\mathcal{C}$ of length $N$ initialized to zero.
For each prime $p \nmid b$, we initialize an update array $\mathcal{P}$
of length $p$ such that the $i^{\mathrm{th}}$ entry of $\mathcal{P}$
is equal to $\log (N_p(b^{-1}(a_0 + i))/p)$.
We then repeatedly add the update array $\mathcal{P}$ into $\mathcal{C}$,
starting with position zero in $\mathcal{C}$ and shifting the starting position
by $p$ with each iteration. Doing this for each prime $p \le B$
tallies the sum $S(t,B)$ into the counter array $\mathcal{C}$
for all $t =a/b$ with $a_0 \le a < a_0 +N$.

By loading $\mathcal{P}$ non-sequentially, we can read the values of  
$\log (N_p(b^{-1}(a_0 + i))/p)$ sequentially from memory,
while requiring only a single inversion modulo $p$ and no additional
multiplications, divisions, or modular reductions.

To avoid the cost of floating point operations, we do not store
$\log(N_p(t)/p)$ as a floating-point number,
but round it to a rational number with fixed denominator~$D$\/
and store the numerator $\lfloor D \log(N_p(t)/p) + \frac12 \rfloor$.
The sieve than tallies these numerators for each~$t$
using integer addition, which is faster than floating-point arithmetic.
The common denominator $D$\/ should be large enough that rounding errors do not
appreciably degrade the score, but small enough that we can keep
a large counter array in the high-speed cache. We found that by taking
$D = 1024$, we were able to fit all of our scores into 16-bit integers.

We further took advantage of a feature of modern processors known as
vector instructions. These are processor level instructions that
can be used to perform the same operation on multiple consecutive elements of
an array simultaneously.  This allowed us to add 16 elements from the
update array $\mathcal{P}$ into the counter array $\mathcal{A}$ at once,
rather than one at a time.

Compared with computing each $S(t,B)$ individually, sieving is extremely fast.
For example, for a fixed value of $b$, we are able to compute
$S(a/b,2^{16})$ for $2^{20}$ values of $a$ in $3.2$ seconds
on a single thread of a hyperthreaded $2.3$ GHz Intel Skylake Xeon processor.
Smaller values of $B$ are even more efficient;
computing $S(a/b,2^{13})$ for $2^{20}$ values of $a$
takes only $0.19$ seconds on the same processor.

The large speed-up offered by this sieve-like technique
is available only in the first step of Nagao's second trick described above:
we can use it to quickly compute $S(t,B_0)$ for all $t$ in the search region,
but not to compute $S(t,B_i)$ for $i \ge 1$ on a restricted set of $t$.
For $i \geq 1$ we must look up individual values of $\log (N_p(t)/p)$.
However, because the sieve-like technique is so efficient,
we can set $B_0$ large enough that computing $S(t,B_0)$ is
the dominant portion of the work --- see Section \ref{sec:choosingBC}.

\section{Choosing Fibrations}
\label{sec:whichsurfaces}

Perhaps the most important ingredient in searching for high-rank
elliptic curves is choosing a good fibration to search on.
We'll describe the factors that guided our choices,
while leaving the specific choices of fibrations to
Section \ref{sec:Z2Z} -- Section \ref{sec:Z7Z}.

In the past, the largest rank elliptic curves having torsion subgroups
$\ZZ/2\ZZ$, $\ZZ/3\ZZ$, and $\ZZ/4\ZZ$ have come from specializations of
K3 surfaces having relatively large rank
($9$ for $\ZZ/2\ZZ$, $5$ for $\ZZ/3\ZZ$, and $4$ for $\ZZ/4\ZZ$).
Our search was no different, focusing on the same families in which
the previous records were found.

By contrast, high-rank K3 surfaces are not known to exist
for the other torsion groups we considered.
The largest known rank of a K3 surface having torsion subgroup
$\ZZ/5\ZZ$ or $\ZZ/6\ZZ$ is $1$,
and the universal elliptic curve having a point of order~$7$
is already a K3 surface, of generic rank zero.
As a result, previous searches have focused on high-degree elliptic surfaces
of larger rank \cite{Lecacheux,DujellaPeralTadic}. 

We initially attempted to do the same for the group $\ZZ/6\ZZ$
using a degree $4$ surface of Kihara having rank $3$ \cite{Kihara}
considered in \cite{DujellaPeralTadic}.
We found that while this surface has a relatively large number of
low-height rank 8 specializations, we could not find any such specializations
of height larger than $\approx 2^{13.5}$.
This suggested that as the height of $t$ grew,
either the number of high-rank specializations in this family decayed rapidly
or our scores quickly became less meaningful.

While \cite{DujellaPeralTadic} considered other degree $4$ surfaces
having \MoW\ group $\ZZ/6\ZZ \times \ZZ^3$,
we concluded that the low-hanging fruit on these had already been discovered,
and that our best hope of finding a rank $9$ curve
having torsion subgroup $\ZZ/6\ZZ$ was
to search on the universal elliptic curve with a point of order~$6$,
which is a rational surface.
We made a similar decision regarding the groups $\ZZ/n\ZZ$ for $n=5$ and $n=7$,
for which the universal elliptic curve over ${\rm X}_1(N)$ is respectively
rational and~K3.

\section{Computing ranks}

After finding a set of values of $t$ such that $S(t,B)$ is sufficiently large,
we are left with the problem of identifying
those that actually have large rank.
We approach this problem in two stages.
First, we use descent methods to obtain an upper bound on the rank.
For those specializations where the upper bound is sufficiently large,
we then search for points on whichever coverings we can efficiently compute.

\subsection{Descent computations}

For half of the families we considered,
the torsion subgroup contains a point of order~$2$,
so we could use Fisher's machinery for computing
rank bounds using $2$-power isogeny Selmer groups,
available in Magma via the command
\texttt{TwoPowerIsogenyDescentRankBound} \cite{FisherHigherDescents}.
For all of the specializations we considered where this upper bound was
at least as large as the previous record in the family,
the upper bound was in fact equal to the rank
(though of course we did not know this until after we searched for points).

For the specializations with torsion subgroup $\ZZ/3\ZZ$,
there is no $2$-isogeny over~$\QQ$,
and a full $2$-descent was out of reach.
This forced us to consider a different approach.

As a first cut, we ran all of the high scorers through
a slightly modified version of Magma's \texttt{ThreeIsogenySelmerGroups}
command to obtain a coarse rank bound. 
While the rank bound coming from 3-descent via isogeny tends to be
reasonably tight for small curves, many of the specializations we considered
had a large number of places of split multiplicative reduction,
which boosted this bound for structural reasons unconnected to rank.
To deal with this, we then used our own implementation of the algorithm for
computing the Cassels-Tate pairing developed by Fisher and van~Beek
\cite{vanBeekThesis,FisherVanBeek}
to compute the \hbox{$3$-Selmer} rank of each specialization for which
the rank bound coming from \hbox{$3$-isogeny} descent was at least $14$.

For the curves with $\ZZ/5\ZZ$ and $\ZZ/7\ZZ$ torsion,
we could not easily compute any isogeny Selmer groups.
However, because the fibration with $\ZZ/5\ZZ$ torsion that we searched
is a rational surface over $\QQ(t)$,
the specializations we were interested in where sufficiently small
that we could compute both the $2$-Selmer group and the
Cassels-Tate pairings for each one.

Na\"{\i}vely, we did not expect the same to be possible for the
$\ZZ/7\ZZ$ fibration we considered, because it is a K3 surface over $\QQ(t)$.
However, we discovered that
while the discriminant of this surface has degree $24$,
the discriminant of its \hbox{$2$-division} field has degree only~$6$.
As a result, although the curves in question were quite large, it was still
possible to perform $2$-descent and the Cassels-Tate pairings on them.

\subsection{Searching for points}

Once we had candidate curves that
our Selmer computations suggested had large rank,
we needed to find enough independent points on them to verify that
they indeed had the expected rank.

Our main method for finding these points was by searching for points on
$2$-coverings of each curve using Magma's built-in functionality.
For most of the groups --- $\ZZ/4\ZZ$, $\ZZ/5\ZZ$, $\ZZ/6\ZZ$, and $\ZZ/7\ZZ$ 
-- we were able to compute the complete $2$-Selmer group for
each of the curves in question.

For the group $\ZZ/2\ZZ$, we did the next best thing,
computing the coverings corresponding to the elements of the
Selmer group of a $2$-isogeny and its dual, and searching on those.

In principle, we could have done something similar with
the $3$-isogeny coverings for the curves having torsion subgroup $\ZZ/3\ZZ$ using Elkies's
lattice-based method of searching for points on cubic curves in $\PP^2$
\cite{ElkiesPtSearch}.
However, due to a memory leak we discovered\footnote{While we discovered the presence of this memory leak, 
we did not attempt to identify its source.} in Magma's implementation of Elkies's method, doing so would have required additional effort.
Instead, we searched the $2$-coverings corresponding to the known points on
each curve coming from the rational points on the surface $\E$,
adding new $2$-coverings to the mix whenever we discovered an additional point.

Somewhat surprisingly, this method worked extremely well.
We suspect that because each of the curves in question has
a large number of points of low height,
we likely would have found them using nearly any method we attempted.

\section{Choosing parameters}
\label{sec:choosingBC}

There is an art to choosing proper values for $B_i$ and $C_i$.
The goal of course is to minimize the total time spent searching,
while not missing any of the top candidates. How to do this is unclear.
We chose our values experimentally, and we suspect that our choices were
far from optimal --- see Section \ref{sec:openQ}. Some tradeoffs however are straightforward.

If $C_0$ is too small, then too many values of $t$ pass the initial cutoff,
so the cost of computing $S(t,B_i)$ for $i \ge 1$ dominates, because
looking up the values of $\log (N_p(t)/p)$ individually
is far more expensive than sieving.  Conversely, if $C_0$ is too large
then we risk eliminating promising values of~$t$.

We compromised by choosing $C_0$ rather aggressively ---
targeting a cutdown on the order of $10^{3}$ ---
but using a large enough value of $B_0$ (between $2^{13}$ and $2^{16}$)
to limit the risk of losing any good candidate $t$.
(Previous searches have tended to take $B < 10^3$,
so this seemed sufficiently conservative.)

The values of $B_i$ for $i \ge 1$ are less important.
We chose the $B_i$ to be successive powers of~$2$ up to $B = 2^{18}$.
We also chose our $C_i$ less aggressively for $i \ge 1$,
since these have a smaller effect on the runtime.

\subsection{Skewed search regions}

For some of the fibrations we considered, the polynomials defining the
non-trivial coefficients of $\E$ were skew in the sense of  \cite{Murphy}.
Very roughly, this means that the higher degree coefficients
tend to have larger magnitude than the smaller ones or vice versa. 

As a result, the average magnitude of the coefficients of an integral model
for $E_t$ on a skewed search region (that is, $t = a/b$ with
$\Max(|a|) = s \Max(|b|)$ for some $s \in \QQ$) will be smaller than
the average magnitude of the coefficients of an integral model for $E_t$
on a square search region having the same size.
While we don't have a firm grasp on how the existence of high-rank
specializations is related to the coefficient size of $E_t$,
it seems sensible to search for smaller curves,
so we skewed our search regions accordingly.

\section{Open questions}
\label{sec:openQ}

Although our search was largely successful,
we are left with some open questions regarding the method of Mestre and Nagao.

\begin{enumerate}[1.]

\item How large a prime bound should we be using
 relative to the search region/degree of the family?

Our experience indicates that the score $S(t,B)$ tends to be
 a poorer indicator of rank as the size of the search region grows,
 and that the rate at which it becomes less useful depends on
 the degree of the surface and on its torsion subgroup.

This is unsurprising, since we expect the convergence rate of
the Euler product for $L_{E_t}(s)$ to depend on the conductor,
which in turn grows roughly as a power of the height $H(t)$
depending on the degree and fiber types of the surface.
(More precisely, the conductor is bounded above by a multiple of
that power of $H(t)$, and for typical~$t$ this is the correct growth order.)
We should therefore expect that we need to allow our prime bound $B$\/
%
%
to grow as a function of $\E$ and $H(t)$
in order for $S(t,B)$ to remain useful.
Is it possible to make this relationship precise?

\item How can we incorporate the Tamagawa factors
  at the places where $E_t$ has bad reduction?

It has been observed that the known curves of high rank tend to have
split multiplicative reduction and large Tamagawa numbers at many small primes.
While the $L$-function includes terms for the bad primes
and these can be incorporated into $S(t,B)$,
these terms don't incorporate the Tamagawa numbers at all
and contribution they would make to $S(t,B)$ would be negative. 

One idea would be to include these primes into the score
via the term $\log \left (\frac{c_p(E_t)(p-1)}{p}\right )$.
However, this seems odd, because for surfaces with an isogeny,
the Tamagawa numbers of $E_t$ and its isogenous curves
will generally not be the same,
and any score that hopes to predict the rank should be isogeny-invariant.

In our searches, we found that including the term
$\log \left (\frac{c(p-1)}{p}\right )$ with $1.2 \le c \le 1.68$ in $S(t,B)$
at each prime of split multiplicative reduction (effectively giving
the specialization a fixed bonus for each prime of multiplicative reduction)
tended to work reasonably well. At the same time, this is clearly a hack,
and it would be nice to understand what is the correct thing to do.

\item How closely should the rank be expected to correlate with $S(t,B)$?

One problem that we struggled with was understanding exactly how
the score $S(t,B)$ should relate to the rank of~$E_t$.  For now,
we are forced to choose our bounds conservatively to avoid missing any
high-rank curves, which results in an increased amount of work,
particularly at the descent steps.

Ideally, we would have a Bayesian score 
$\Prob(E_t \text{ has rank at least } r \mid S(t,B) > C)$
that would let us set the bounds $B_i$ and $C_i$ optimally,
and inform our decision about how many curves to apply descent methods to.
(The use of a Bayesian score was suggested to us by Joel Rosenberg.) 
Such a score would also let us estimate the likelihood that
we missed a curve of high rank.

\end{enumerate}

\section{Main results}
\label{sec:MainResults}

We obtained new rank records for elliptic curves with torsion subgroups
$\ZZ/2\ZZ$,  $\ZZ/3\ZZ$, $\ZZ/4\ZZ$,$\ZZ/6\ZZ$, and $\ZZ/7\ZZ$.
We also carried out an extremely large search for a record-breaking curve
having torsion subgroup $\ZZ/5\ZZ$, but with no success.

The current and previous records (as given by \cite{Dujella})
for each of these torsion subgroups are given in Table \ref{tab:recordtable}.
The next sections describe in greater detail the searches we carried out
in pursuit of each of these records.

\begin{table}[h]
\renewcommand{\arraystretch}{1.25}
\begin{tabular}{c|c|c}
Torsion Subgroup & Previous Record & Current Record \\
\hline
$\ZZ/2\ZZ$ & 19 & 20 \\ 
$\ZZ/3\ZZ$ & 14 & 15 \\ 
$\ZZ/4\ZZ$ & 12 & 13 \\ 
$\ZZ/5\ZZ$ & 8 & 8 \\ 
$\ZZ/6\ZZ$ & 8 & 9  \\ 
$\ZZ/7\ZZ$ & 5 & 6 \\ 
\end{tabular}
\vspace{3pt}
\caption{Rank records for various torsion subgroups}
\label{tab:recordtable}
\end{table}



\section{Curves with torsion subgroup $\ZZ/2\ZZ$}
\label{sec:Z2Z}

For torsion groups $T = \ZZ/2\ZZ$, $\ZZ/3\ZZ$, $\ZZ/4\ZZ$
we proceeded as in~\cite{Elkies3Lectures},
computing an elliptic fibration $\E(\QQ_t)$ of a K3 surface~$X$\/
whose N\'{e}ron-Severi group $\NS(X)$ is defined over~$\QQ$ and
has high rank and large discriminant.
For $T=\ZZ/3\ZZ$ and $T=\ZZ/4\ZZ$ we used the surface with
$\NS(X)$ of rank~$20$ and discriminant $-163$.
But for $T=\ZZ/2\ZZ$ this discriminant is not large enough;
it turns out \cite{ElkiesCensus} that the highest rank attained by
an elliptic fibration of~$X$\/ with a \hbox{$2$-torsion point} is~$8$.
Instead we use $X$\/ with $\NS(X)$ of rank~$19$ but larger discriminant,
which can attain \MoW\ rank~$9$.

Such $X$\/ are parametrized by elliptic or Shimura modular curves,
call them~$C$, of level $\frac12\left|\disc\,\NS(X)\right|$.  When 
$\left|\disc\,\NS(X)\right|$ is large enough to allow
\MoW\ rank~$9$, the curve $C$\/ usually has genus at least~$2$,
with few if any rational points (other than cusps and CM points,
at which $X$\/ or the elliptic fibration degenerates).
In~\cite[p.8--9]{Elkies3Lectures} Elkies reports using the
sporadic rational point on the {genus-$2$} curve ${\rm X}_0(191)/w$
to find such an~$X$.  A few years later he found a genus-zero
Shimura curve of level $230$ that could be used instead,
giving a family of elliptic surfaces with \MoW\ group
$\ZZ/2\ZZ \times \ZZ^9$.  Here $C = {\mathcal X} / w_{230}$,
with $\mathcal X$\/ associated to the congruence subgroup $\Gamma_0(23)$
of the quaternion algebra ramified at $\{2,5\}$.
The family of surfaces with their elliptic fibrations
was computed as in \cite{ElkiesShimuraK3,EK}.
The elliptic fibration is of the form
$\E_u / \QQ(t) : y^2 = x^3 + 2Ax^2 + Bx$, where
\begin{multline}
\label{eq:Eu_A}
A = (u^8-18u^6+163u^4-1152u^2+4096) t^4
   + (3u^7-35u^5-120u^3+1536u) t^3 \\
   + (u^8-13u^6+32u^4-152u^2+1536) t^2 
   + (u^7+3u^5-156u^3+672u) t \\
   + (3u^6-33u^4+112u^2-80),
\end{multline}
and $B = \prod_{i=1}^8 B_i(t,u)$ where
\begin{multline}
\label{eq:Eu_B}
B_1(t,u) = (u^2+u-8)t + (-u+2), \quad
B_3(t,u) = (u^2-u-8)t + (u^2+u-10), \\
B_5(t,u) = (u^2-7u+8)t + (-u^2+u+2), \quad
B_7(t,u) = (u^2+5u+8)t + (u^2+3u+2),
\end{multline}
and $B_i(t,u) = -B_{i-1}(-t,-u)$ for $i=2,4,6,8$.
Thus $\E_u \cong \E_{-u}$.  If $5-u^2$ is a square, and
$u \neq \pm 1,\pm 2$ (to exclude CM points), then
$\E_u$ has \MoW\ group $\ZZ/2\ZZ \times \ZZ^9$ over $\QQ(t)$.
Generators are exhibited in Appendix~\ref{app:pointsonShimura}.

We searched for high-rank specializations of $\E_u$ for several values of $u$.

For $u = 2/5$, we searched the region
$t = a/b$ with $0 < a < 2^{21}$ and $-2^{23} < b < 2^{23}$,
finding 17 curves of rank 19, including the previous record-holding curve
of Elkies that appears in \cite{Dujella}, which occurs at $t = 11860/97527$.

For $u = 11/5$, we first applied the linear fractional transformation
$t \mapsto \frac{2-t}{t-6}$ to $\E_u$ and then searched the region
$t = a/b$ with $0 < a < 3 \cdot 2^{21}$ and $-2^{21} < b < 2^{21}$.
We found one specialization of rank 20 at $t = -68559/32629$
($t = -721141 / 2026305$ on the orginal model of $\E_u$),
as well as another 20 specializations of rank 19,
including one at $t = 100782/104143$
($t = -26876/131019$ on the orginal model of $\E_u$)
with smaller discriminant than the rank 19 curve of Elkies
appearing in~\cite{Dujella}.

Minimal models and $x$-coordinates of a set of generators
for the rank 20 specialization and the smallest conductor rank 19
specialization appear in Appendix \ref{app:xcoords2}.
We note that this curve of rank 20 is the elliptic curve of largest rank
for which the rank is known unconditionally.

We also searched regions of size roughly $2^{44}$
on each of the fibrations coming from $u = 2/13$ and $u = 22/13$,
but did not find any specializations of rank greater than 18.

\section{Curves with torsion subgroup $\ZZ/3\ZZ$}

We searched the 13 K3 surfaces having \MoW\ group
$\ZZ/3\ZZ \times \ZZ^5$ families that occur as elliptic fibrations on
the singular K3 surface of discriminant $-163$
and which will appear in \cite{ElkiesCensus}.

We searched an appropriately skewed region of size $2^{43}$
on each of the 13 fibration, finding 34 specializations of rank 14
(at least one on 11 of the 13 fibrations)
as well as a single specialization of rank 15, given by
$$
E:y^2 + 490738465519xy - 432802729180188878035670522423557875y = x^3.
$$
Among the specializations having rank $14$, the one with smallest
conductor and discriminant is given by
$$
y^2 + 6244332976xy + 2204421250641922174556630375y = x^3,
$$
which has smaller conductor and discriminant than the previously known
curve of rank 14 appearing in~\cite{Dujella}.
The $x$-coordinates of a set of generators for each of these curves
is given in given in Appendix~\ref{app:xcoords3}.

\section{Curves with torsion subgroup $\ZZ/4\ZZ$}

We searched a pair of families each having \MoW\ group
$\ZZ^4 \times \ZZ/4\ZZ$, both of which are
elliptic fibrations of the singular K3 surface of discriminant $-163$.
The first, given by the equation
\begin{multline}\E_1: y^2 + (8t-1)(32t+ 7)xy + 8(8t-1)(32t+7)(t+1)(15t-8)(31t-7)y \\ = x^3 +  8(t+1)(15t-8)(31t-7)x^2,\end{multline}
appears (with a typo) in \cite{Elkies3Lectures};
the second, given by the equation
\begin{multline}
\E_2: y^2 -8 (80t+9) - 16(80t+9)(t-2)(2t-1)(18t-1)(2t-81)y \\ = x^3 +  2 (t-2)(2t-1)(18t-1)(2t-81)x^2,
\end{multline}
will appear in \cite{ElkiesCensus}.

The previous rank record for torsion group $\ZZ/4\ZZ$ was $12$,
attained by two curves in the family $\E_1$, found by
Elkies in 2006 ($t = 18745/6321)$
and Dujella-Peral in 2014 ($t = -13083/72895)$.
Searching up to height $2^{22}$ on $\E_1$,
we found three rank 13 specializations
at $t=-1086829/638219$, $t=-2856967/190447$, and $t = 973215/3135431$,
as well as 76 rank 12 specializations.
Among the rank 12 specializations, the one with smallest conductor occurs at
$t=-447577/2601952$ ($\mathfrak{f}_{E_t}~\approx~2^{153.41}$)
and the one with smallest discriminant occurs at
$t=83497/251378$ ($|\Delta_{E_t}| \approx 2^{392.96}$).
Respectively, these have smaller conductor and discriminant than
the previously known rank 12 curves.

We searched up to height $2^{22}$ on $\E_2$
and were unable to find any specializations of rank 13,
though we did find 32 having rank 12. 

\section{Curves with torsion subgroup $\ZZ/5\ZZ$}

As noted in Section \ref{sec:whichsurfaces},
for the group $\ZZ/5\ZZ$, we chose to search for good specializations
on the universal elliptic curve having a point of order~$5$,
which is a rational elliptic surface. A model for this surface is given by 
$$
y^2 + (t + 1)xy + ty = x^3 + tx^2
$$
and we searched for $t$ up to height $2^{28}$ on this surface.
Unfortunately, our search did not yield any any specializations
having rank larger than the current record of 8. 

We did however find 151 rank 8 specializations, three of which were
previously known. The curve we found with smallest conductor appears at
$t = 1809535/5292661$ ($\mathfrak{f}_{E_t} \approx 2^{85.86}$
and the curve we found with smallest discriminant appears at
$t = 5167107/723695$ ($|\Delta_{E_t}| \approx 2^{254.77}$).
Each of these has both smaller conductor and discriminant than
all of the previously known rank 8 curves.

Minimal models and $x$-coordinates for a set of generators for
each of these curves are given in 
given in Appendix \ref{app:xcoords5}.
 
\section{Curves with torsion subgroup $\ZZ/6\ZZ$}
As was the case for $\ZZ/5\ZZ$, we chose to search for good specializations
on the universal elliptic curve having a point of order~$6$,
which is a rational elliptic surface. A model for this surface is given by
$$
y^2 + txy + (t+2)y = x^3,
$$
with torsion points of order $2,3,6$ at
$(x,y) = (-1,-1), (0,0), (t+2,t+2)$ respectively.

We searched for good specializations of this model in the region 
$t = a/b$ with $0 < a < 2^{25}$ and $-2^{26} < b < 2^{26}$.
In this case, the skewed search region was a fortuitous accident,
rather than a deliberate choice.
We found a single rank~$9$ curve at $t = -22029701/37178488$
as well as $71$ rank~$8$ specializations,
all but one of which appear to be previously unknown.
The rank~$8$ curve with the smallest conductor and smallest discriminant
appears at $t = 6308333/1000939$ ($\mathfrak{f}_{E_t} \approx 2^{81.96}$ and
$|\Delta_{E_t}| \approx 2^{253.07}$). Its $2$-isogenous curve that appears at
$t = -24627934/8310211$ shares the same conductor, but has larger discriminant.

Minimal models and $x$-coordinates of a set of generators
for the rank~$9$ specialization
and the smallest conductor/discriminant rank~$8$ specialization
appear in Appendix~\ref{app:xcoords6}.

\begin{remark}
In retrospect, we could have taken advantage of the involution
$w_2: t \mapsto -(2t+12)/(t+2)$, for which $E_{w_2 t}$ is the curve
$E_t'$ which is \hbox{$2$-isogenous} with $E_t$,
and thus also has torsion subgroup $\ZZ/6\ZZ$.
This would let us restrict our search area to $-4 < t < 2$.
In partial compensation, we could compare the scores of $t$ and $w_2(t)$
to corroborate that we are computing these scores correctly.

\end{remark}

\section{Curves with torsion subgroup $\ZZ/7\ZZ$}
\label{sec:Z7Z}

As noted in Section \ref{sec:whichsurfaces},
for the group $\ZZ/7\ZZ$, we chose to search for good specializations of
the universal elliptic curve having a point of order~$7$,
Unlike the groups $\ZZ/5\ZZ$ and $\ZZ/6\ZZ$,
the universal elliptic curve having a point of order~$7$
is a K3 surface rather than a rational one. 

A model for this curve is given by 
$$
y^2 + (-t^2 + t + 1)xy + (-t^3 + t^2)y = x^3 + (-t^3 + t^2)x^2
$$ 
(see e.g.\ \cite[p.195]{Tate:AEC}).
We searched up to height $2^{20}$ on this model and found
a single specialization of rank~$6$ at $t = -748328/820369$.
A minimal model and the set of $x$-coordinates of a set of generators of
this specialization are given in Appendix \ref{app:xcoords7}.

\begin{remark}
We also attempted to search for record breaking curves in the K3 families
having torsion subgroups $\ZZ/8\ZZ$ and $\ZZ/2\ZZ \times \ZZ/6\ZZ$
without any luck. We suspect that our success here is due to
$\ZZ/7\ZZ$ having been the lone torsion subgroup for which
the generic elliptic curve is K3 but for which there was
no previously known curve having rank~$6$.
\end{remark}

\appendix

\section{Points on $\E_u/\QQ(t)$}
\label{app:pointsonShimura}

Recall that in (\ref{eq:Eu_A},\ref{eq:Eu_B})
we exhibit $A$ and $B_1,\ldots,B_8$ in $\QQ[t,u]$ such that
$\E_u/\QQ(t)$ has Weierstrass equation $y^2 = x^3 + 2Ax^2 + Bx$
where $B = \prod_{i=1}^8 B_i$.  The minimal height of
a non-torsion section is~$2$, attained by $70$ pairs $(x,\pm y)$
with $x,y \in \QQ(u,\sqrt{5-u^2})[t]$.
We find that $58$ of the $70$ pairs have $x,y \in \QQ(u)[t]$; these generate a
\MoW\ subgroup of rank~$8$.  One simple choice of generators of this subgroup
consists of points with \hbox{$x$-coordinates}
\begin{multline}
\label{eq:x1-x8}
  -\!B_1 B_2 B_3 B_6,\ -B_1 B_2 B_4 B_5,\ 4 B_1 B_2 B_5 B_6,\ B_1 B_3 B_4 B_6,
\\
  -\!B_1 B_3 B_4 B_7,\ B_1 B_3 B_4 B_8,\ B_1 B_3 B_5 B_6,\ -B_1 B_5 B_6 B_7.
\end{multline}
Extending $\QQ(u)$ by $\sqrt{5-u^2}$ yields $\QQ(m)$ where
$m$ is a rational coordinate on the parametrizing Shimura curve, with
\begin{equation}
\label{eq:m,u}
u = 2 \, \frac{m^2-m-1}{m^2+1},\quad
(5-u^2)^{1/2} = \pm\frac{m^2+4m-1}{m^2+1};
\end{equation}
then adding $-(m-1)^2 B_1 B_2 B_3 B_8$ to the list (\ref{eq:x1-x8})
gives \hbox{$x$-coordinates} of $9$ \MoW\ generators modulo torsion.
The Gram matrix of canonical height pairings is
\begin{equation}
\label{eq:Gram9}
\frac12 \left[
\begin{array}{rrrrrrrrr}
 4 & 0 & 1 &-1 & 0 & 2 &-1 & 0 & 1 \cr
 0 & 4 &-1 &-2 & 0 & 2 &-2 & 0 & 0 \cr
 1 &-1 & 4 & 0 &-1 & 1 &-1 & 1 & 2 \cr
-1 &-2 & 0 & 4 &-1 &-1 & 1 & 0 & 0 \cr
 0 & 0 &-1 &-1 & 4 & 1 & 0 &-2 & 0 \cr
 2 & 2 & 1 &-1 & 1 & 4 &-2 &-1 & 1 \cr
-1 &-2 &-1 & 1 & 0 &-2 & 4 & 1 & 0 \cr
 0 & 0 & 1 & 0 &-2 &-1 & 1 & 4 & 1 \cr
 1 & 0 & 2 & 0 & 0 & 1 & 0 & 1 & 4
\end{array}
\right],
\end{equation}
with determinant $115/16$.

\section{Models for record breaking curves} 
\label{app:modelsandpoints}

\subsection{Overview}
This section gives minimal integral models for each of
the record breaking curves we discovered,
along with the $x$-coordinates of a set of generators for
the torsion-free part of each of them.

By common convention we use a vector $(a_1,a_2,a_3,a_4,a_6)$
to mean the extended Weierstrass model
$$
y^2 + a_1 x y + a_3 x = x^3 + a_2 x + a_4 x + a_6
$$
whose coefficients are the vector's entries.
We usually depart from another common convention that chooses the model with
$a_1,a_3 \in \{0,1\}$ and $a_2 \in \{-1,0,1\}$.
Such models have the advantage of being unique,
but for curves with nontrivial torsion there may be
one or more other choices that put a torsion point at $(x,y)=(0,0)$
and have a coefficient vector with noticeably fewer digits 
(for starters $a_6=0$ if $(0,0)$ is on the curve).

When possible we give a generating set of $E(\QQ) \bmod E(\QQ)_{\rm tors}$
consisting of integral points of small height.  For most of our curves
there are plenty of such points to choose from, even though there can be
other curves with the same torsion group and somewhat lower rank
that have even more integral points.

\subsection{$\ZZ/2\ZZ$}
\label{app:xcoords2}

A minimal model for the rank $20$ curve having  $\ZZ/2\ZZ$ torsion
has coefficients
\begin{multline*}
(1, -1, 1, -244537673336319601463803487168961769270757573821859853707, \\
961710182053183034546222979258806817743270682028964434238957830989898438151121499931).
\end{multline*}
Here we reluctantly give a model with small $a_1,a_2,a_3$ and huge $a_4,a_6$,
because the torsion point has $x = -69288588686111702678625616725/4$
and thus cannot be put at the origin on a minimal model.\footnote{
  The coefficients
  ${\scriptstyle
  (2,-207865766058335108035876850179,0,
   10490122792958386322093670444427223877319227761081795217921,0)
  }
  $
  give a model with smaller coefficients that puts the torsion point at $(0,0)$
  but is not minimal at~$2$.
  }

One choice of $20$ points that generate its \MoW\ group modulo torsion
has \hbox{$x$-coordinates}
$$
\begin{array}{rr}
  -5976635286513806621064126789, &     595416388787490259443766591, \cr
   2434562872293108275107029075, &    3513074027344435171140978981, \cr
    399682145249051758133327419, &  -10714754038296881855524018251, \cr
 -16034220456847626275437501599, &    1185828672355214392425799131, \cr
 -11190697582885409770718510409, &    2634316446310680332042122261, \cr
  64222149978369055569434725591, &   23945425437351916471937562579, \cr
  13094114400583295432756346651, & 2689776334541089917424552236511, \cr
  -2627014038979941829331861469, &  113605800622499112413124359631, \cr
  -7364938748841807757773625709, &  -14298222927159284914180072349, \cr
 785686589410787916270883192839, &   -2250170491079839258934900709.
\end{array}
$$
Here and later we list generators in increasing order by canonical height.

A minimal model for the rank $19$ curve with $\ZZ/2\ZZ$ torsion
having smallest known conductor has coefficients
\begin{multline*}
(1, 4040549489437705068551042, 0,
  39096673111815206065773237234587256582331296000, 0)
\end{multline*}

One choice of $19$ points that generate its \MoW\ group modulo torsion
has \hbox{$x$-coordinates}

$$
\begin{array}{rr}
-3613294426098135199878600, &   284077053735716552925900, \cr
  -69786343891815820666800, &     6409078899434870587500, \cr
 4711243262341394854929360, &  -200862034480295787990300, \cr
   49746704013683926431600, &     1283007628272047952000, \cr
  601243680664306184613420, &  1681679070386109358006014, \cr
 -178674347439204200162150, &  -140058466067600728971180, \cr
    4490592251930741573760, & -1245418009246864352006250, \cr
  239435938047242410050720, & -2615926042511102882808000, \cr
-3662820474106418641536000, &  308679854892675472378120, \cr
  -12130119373140047385600.
\end{array}
$$

\subsection{$\ZZ/3\ZZ$}
\label{app:xcoords3}

The rank $15$ elliptic curve with coefficient vector
$$
(490738465519, 0, -432802729180188878035670522423557875, 0, 0)
$$
has a $3$-torsion point at $(x,y)=(0,0)$.
One choice of $15$ points that generate its \MoW\ group modulo torsion
has \hbox{$x$-coordinates}

\vspace*{1ex}

\centerline{
$
\begin{array}{rrr}
 414082294873186000299147, & -461076037958619691375950, &   136016697778663191410466, \cr
 579811074194569447550775, & 4156065765459153070875350, &  -379256436856490083222605, \cr
-480257266200757201099125, &  626879349686994759271350, &   319402198167922579675875, \cr
9987762741068630814895872, & 1025559076978453798187316, & 17710047123788181654048375, \cr
 236426830570889446065942, & -162860681446721622110565, &  1093411474853808475876875.
\end{array}
$
}

The rank $14$ elliptic curve with coefficient vector
$$
(6244332976, 0, -2204421250641922174556630375, 0, 0)
$$
has a $3$-torsion point at $(x,y)=(0,0)$.
One choice of $14$ points that generate its \MoW\ group modulo torsion
has \hbox{$x$-coordinates}

\vspace*{1ex}

\centerline{
$
\begin{array}{rrrr}
2907919170263662, & -65199074165293250, & 71604990115331040,  & 77567806466944000, 
\cr 
108999498650081840, & 169617569990697350, & -171009947870163008, & -204167066230390100,
\cr 
-240427032442334750, & 243676691791782250, & -256142889038646510, & -276580713950955750, 
\cr 
368313341140417750, & -449841531945448000.
\end{array}
$}

\subsection{$\ZZ/4\ZZ$}
\label{app:xcoords4}
The first rank $13$ curve with $\ZZ/4\ZZ$ torsion has a minimal model
with coefficient vector
\begin{multline*}
(282887999996745, -1871148179781457712818452480, \\
  -529325366275926422138597740307015937177600, 0, 0)
\end{multline*}
and a $4$-torsion point at $(x,y)=(0,0)$.
One choice of $13$ points that generate its \MoW\ group modulo torsion
has \hbox{$x$-coordinates}

$$
\begin{array}{rr}
 37563104221873287230436120000, &    1241851783771179145432296000, \cr
  1992140999686088390294877150, &   30921042737991542683359263880, \cr
-21195532433936174709304166400, &   -1464098167733086800531916800, \cr
  1670745991840921221771294750, &    1252355926117744178967180450, \cr
 -1960920553671074388872220170, & 1375293185347275499663130572800, \cr
  2549902537861429590505036800, & 3272919221738028252106303872714, \cr
102225511700163143939329914880.
\end{array}
$$

The second rank $13$ curve with $\ZZ/4\ZZ$ torsion has a minimal model
with coefficient vector

\begin{multline*}
(230691818102905, -200100346570723590045845120, \\
-46161512753421616727023025112895852073600, 0, 0)
\end{multline*}
and a $4$-torsion point at $(x,y)=(0,0)$.
One choice of $13$ points that generate its \MoW\ group modulo torsion
has \hbox{$x$-coordinates}

$$
\begin{array}{rr}
   190412869629748629206788500, & -11655521125151390350616252280, \cr
-10482658909728296079200226100, &   -205253870232797421109008000, \cr
   193230556828647163522857600, &   2390337099874364874239977850, \cr
-10561431236301791011714683300, &  -1195165694989063921020955200, \cr
   876665740401972718169616600, &    -99112055810721390011710344, \cr
   -65566000913948267196883584, &    166949951644450209072942720, \cr
   -26328612670314620364001050.
\end{array}
$$

The third rank $13$ curve with $\ZZ/4\ZZ$ torsion has a minimal model
with coefficient vector
\begin{multline*}
(246888014319233, -8884285566590219865500325632, \\ -2193423622180481268696018169961040300480256, 0, 0 )
\end{multline*}
and a $4$-torsion point at $(x,y)=(0,0)$.
One choice of $13$ points that generate its \MoW\ group modulo torsion
has \hbox{$x$-coordinates}

$$
\begin{array}{rr}
 -968516084234641058709370232, & -1333726837303108113451614080, \cr
 1792794868671671366043266816, &  2362595876319902581142656768, \cr
-2746004168634841009972934984, &  3469325866293712913010729024, \cr
 3644805279133239447459855232, &  4449372053406414078540323280, \cr
-4537829698895530474950049368, &  5156996081584183666047796032, \cr
 5789474008645490085082165824, &  5912795841516183863849831680, \cr
10555676267250916670215460568.
\end{array}
$$

\subsection{$\ZZ/5\ZZ$}
\label{app:xcoords5}

The rank~$8$ curve with $\ZZ/5\ZZ$ torsion
having smallest known conductor has a minimal model
with coefficient vector
$$
(7102196, 9577255322635, 50689165733152681735, 0, 0).
$$
The torsion group is generated by $(x,y) = (0,0)$.
One choice of $8$ points that generate the \MoW\ group modulo torsion
has \hbox{$x$-coordinates}

\begin{equation}
\begin{array}{rrrr}
 -11217531799903, & -10836503720185, & -4357099419673, &  1401549559410, \cr
 256939125827615, & -10247328030940, & -6060818514894, & -6697297034428
\end{array}
\end{equation}

The rank~$8$ curve with $\ZZ/5\ZZ$ torsion
having smallest known discriminant has a minimal model
with coefficient vector
$$
(5890802, 3739409500365, 2706191958366648675, 0, 0).
$$
The torsion group is generated by $(x,y) = (0,0)$.
One choice of $8$ points that generate the \MoW\ group modulo torsion
has \hbox{$x$-coordinates}

\begin{equation}
\begin{array}{rrrr}
 -21207376737, & 37660080920, & -89104376475, & 100531079550, \cr
 117291419735, & -120660570135, & 148808336985, &-214614453600
\end{array}
\end{equation}

\subsection{$\ZZ/6\ZZ$}
\label{app:xcoords6}

The rank~$9$ curve with $\ZZ/6\ZZ$ torsion
has a minimal model with coefficient vector
$$
(-22029701, 0, 72328851024410157777600, 0, 0).
$$
The torsion group is generated by
$(x,y) = (1945448965660200, 72328851024410157777600)$;
multiplying this point by $2$ yields the \hbox{$3$-torsion} point $(0,0)$.
One choice of $9$ points that generate the \MoW\ group modulo torsion
has \hbox{$x$-coordinates}
\begin{equation}
\begin{array}{rrr}
  749629491053742, &   6092756193428190, & -1380249411088240, \cr
-1067429532233440, & 174532909579773030, &   949536320242950, \cr
 1079473135677300, &  24157188371048640, &  3112751229126000.
\end{array}
\end{equation}
 
The rank~$8$ curve with $\ZZ/6\ZZ$ torsion and smallest known conductor and discriminant
has a minimal model with coefficient vector
$$
(6308333, 0, 8325824903545553131, 0, 0).
$$
The torsion group is generated by
$(x,y) = (8318014288129 : 8325824903545553131)$;
multiplying this point by $2$ yields the \hbox{$3$-torsion} point $(0,0)$.
One choice of $8$ points that generate the \MoW\ group modulo torsion
has \hbox{$x$-coordinates}
\begin{equation}
\begin{array}{rrrr}
-204062889121, & 211687889245, & -403788801990, & -410295468023, \cr
-733395115518, & -823562706096, & -859172099915, & -2828410292799.
\end{array}
\end{equation}

\subsection{$\ZZ/7\ZZ$}
\label{app:xcoords7}

The rank~$6$ curve with $\ZZ/7\ZZ$ torsion
has a minimal model with coefficient vector
$$
(-500894592455, 720663120331059917723712,
 485010096730715360294683087532269632, 0, 0).
$$
The torsion group is generated by $(x,y) = (0,0)$.

One choice of $6$ points that generate the \MoW\ group modulo torsion
has \hbox{$x$-coordinates}
\begin{equation}
\begin{array}{rl}
 -863240219455759708343872, & 147841500613888155442368, \cr
 -655405721270483784258504, & 227328163133810400709740, \cr
17758591139156733971281176, & 4457894404162347392127765558505920 / 79519^2.
\end{array}
\end{equation}
The large final generator is inevitable:
the first five generators have canonical heights between $15.434$ and $19.431$,
but the last generator must have height at least $42.058$
(we have made the minimal choice, and with the smallest possible denominator
among its seven torsion translates).

\end{document}